\documentclass[graybox]{svmult}

\usepackage{mathptmx}       
\usepackage{helvet}         
\usepackage{courier}        
\usepackage{type1cm}        
\usepackage{makeidx}         
\usepackage{graphicx}        
\usepackage{multicol}        
\usepackage[bottom]{footmisc}

\usepackage{psfrag}
\usepackage{amsfonts,amsmath,epsfig,amssymb,color,mathrsfs}
\newtheorem{teor}{Theorem}[section]
\newtheorem{defi}[teor]{Definition}

\makeindex             

\begin{document}

%
%
%
%

\title*{Conditional stability for backward parabolic equations with Osgood coefficients}

\author{Daniele Casagrande, Daniele Del Santo and Martino Prizzi}

\institute{Daniele Casagrande \at Universit\`a degli Studi di Udine, Via delle Scienze, 206 - 33100 Udine, Italy, \email{daniele.casagrande@uniud.it}
\and Daniele Del Santo \at Universit\`a degli Studi di Trieste, Via A. Valerio, 10 - 34100 Trieste, Italy \email{delsanto@units.it}
\and Martino Prizzi \at Universit\`a degli Studi di Trieste, Via A. Valerio, 10 - 34100 Trieste, Italy \email{mprizzi@units.it}}

\maketitle

\abstract{
The interest of the scientific community for the existence, uniqueness and stability of solutions to PDE's is testified by the numerous works available in the literature. In particular, in some recent publications on the subject~\cite{NONLINAL,MATAN} an inequality guaranteeing stability is shown to hold provided that the coefficients of the principal part of the differential operator are Log-Lipschitz continuous. Herein this result is improved along two directions. First, we describe how to construct an operator, whose coefficients in the principal part are not Log-Lipschitz continuous, for which the above mentioned inequality does not hold. Second, we show that the stability of the solution is guaranteed, in a suitable functional space, if the coefficients of the principal part are Osgood continuous.}

\section{Introduction}
Backward parabolic equations are known to generate ill-posed (in the sense of Hadamard~\cite{Had_53,Had_64}) Cauchy problems. Due to the smoothing effects of the parabolic operator, in fact, it is not possible, in general, to guarantee existence of the solution for initial data in any reasonable function space. In addition, even when solutions possibly exist, uniqueness does not hold without additional assumptions on the operator. Nevertheless, also for problems which are not well-posed the study of the conditional stability of the solution -- the surrogate of the notion of ``continuous dependence'' when existence of a solution is not guaranteed -- is interesting . Such kind of study can be performed by resorting to the notion of \emph{well behaving} introduced by John~\cite{John}: a problem is \emph{well-behaved} if ``only a fixed percentage of the significant digits need be lost in determining the solution from the data''. In other words, a problem is well behaved if its solutions in a space $\mathcal{H}$ depend continuously on the data belonging to a space $\mathcal{K}$, provided they satisfy a prescribed bound in a space $\mathcal{H}^\prime$ (possibly different from $\mathcal{H}$).
In this paper we give a contribution to the study of the (\emph{well}) behaviour of the Cauchy problem associated with a backward parabolic operator. In particular, we consider the operator $\mathcal{L}$ defined, on the strip $[0,T]\times \mathbb{R}^n$, by
\begin{equation}\label{eq_L}
\mathcal{L}u=\partial_t u +\sum_{i,j=1}^n \partial_{x_i}\left( a_{i,j}(t,x)\partial_{x_j}u \right)+\sum_{j=1}^n b_j(t,x)\partial_{x_j}u+c(t,x)u\,,
\end{equation}
where all the coefficients are bounded. We suppose that $a_{i,j}(t,x)=a_{j,i}(t,x)$ for all $i,j=1,\ldots,n$ and for all $(t,x)\in [0,T]\times\mathbb{R}^n$. We also suppose that $\mathcal{L}$ is backward parabolic, i.e. there exists $k_A\in ]0,1[$ such that, for all $(t,x,\xi)\in[0,T]\times \mathbb{R}^n\times \mathbb{R}^n$,
\begin{equation}
k_A\vert \xi\vert^2\le \sum_{i,j=1}^n a_{i,j}(t,x)\xi_i\xi_j\le k_A^{-1}\vert\xi\vert^2\,.
\end{equation}

We show that if the coefficients of the principal part of $\mathcal{L}$ are at least Osgood regular, then there exists a function space in which the associated Cauchy problem 
\begin{equation}\label{eq_probl_cauchy}
\left\{\begin{array}{ll}
\mathcal{L}u=f\,,\qquad & \textnormal{ in }(0,T)\times \mathbb{R}^n\,,\\
u\vert_{t=0}=u_0\,,\qquad & \textnormal{ in }\mathbb{R}^n\,,
\end{array}\right.
\end{equation}
has a stability property.

To collocate the new result in the framework of the existing literature, we first recall the contents of some interesting publications on the subject which show that, as one could expect, the function space in which the stability property holds is related to the degree of regularity of the coefficients of $\mathcal{L}$. Weaker requirements on the regularity of the coefficients must be balanced, for the stability property to hold, by stronger \emph{a priori} requirements on the regularity of the solution, hence stability holds in a smaller function space.

The overview on available works helps to lead the reader to the new result, claimed in the final part of the paper, concerning operators with Osgood-continuous coefficients. This kind of regularity is critical since it is the minimum required regularity that guarantees uniqueness of the solution and can therefore be considered as a sort of lower limit. The complete proof of the claim is rather cumbersome and is not reported here; instead, we provide a detailed discussion on the fact that, although the core reasoning is based on the theoretical scheme followed to achieve previous results~\cite{MATAN}, the modifications needed to obtain an analogous proof in the case of Osgood coefficients are by no means trivial.

\section{Uniqueness and non-uniqueness results}
We begin by recalling some results on the uniqueness and non-uniqueness of the solution of the problem (\ref{eq_probl_cauchy}). Consider the space
\begin{equation}\label{eq_insieme}
\mathcal{H}_0=C([0,T],L^2)\cap C([0,T[,H^1)\cap C^1([0,T[,L^2)\,.
\end{equation}
One of the first results concerning uniqueness is due to Lions and Malgrange~\cite{Lio_Mal}. They achieve a uniqueness result for an equation associated to a sesquilinear operator defined in a Hilbert space. With respect to the space (\ref{eq_insieme}), this result can be read as follows.

\begin{teor}\label{teo_Lio_Mal}
If the coefficients of the principal part of $\mathcal{L}$ are Lipschitz continuous with respect to $t$ and $x$, if $u\in\mathcal{H}_0$ and if $u_0=0$, then $\mathcal{L}u=0$ implies $u\equiv 0$.$\hfill\square$
\end{teor}
The Lipschitz continuity of the coefficients is a crucial requirement for the claim, as shown some years later by Pli\'s~\cite{Pli} in the following theorem.

\begin{teor}
There exist $u$, $b_1$, $b_2$ and $c\in C^\infty(\mathbb{R}^3)$, bounded with bounded derivatives and periodic in the space variables and there exist $l:[0,T]\to\mathbb{R}$, H\"older-continuous of order $\delta$ for all $\delta<1$ but not Lipschitz-continuous, such that $1/2\le l(t)\le 3/2$ and the support of the solution $u$ of the Cauchy problem
\begin{equation}\label{eq_cauchy_plis}
\left\{\begin{array}{ll}
\partial_t^2u(t,x_1,x_2)+\partial_{x_1}^2u(t,x_1,x_2)+l(t)\partial_{x_2}^2u(t,x_1,x_2)+\\
\qquad +b_1(t,x_1,x_2)\partial_{x_1}u(t,x_1,x_2)+b_2(t,x_1,x_2)\partial_{x_2}u(t,x_1,x_2)+\\
\qquad\qquad +c(t,x_1,x_2)u(t,x_1,x_2)=0 & \textnormal{in }\mathbb{R}^3\,,\\[2mm]
\left.u(t,x_1,x_2)\right\vert_{t=0}=0 & \textnormal{in }\mathbb{R}
\end{array}\right.\!\!\!
\end{equation}
is the set $\mathbb{R}\times\mathbb{R}\times\{t\ge 0\}$.$\hfill\square$
\end{teor}
Note that the differential operator in (\ref{eq_cauchy_plis}) is elliptic. However, the same idea developed by Pli\'s to prove the claim can be exploited to obtain a counterexample for the backward parabolic operator
$$
\mathcal{L}_P=\partial_t+\partial_{x_1}^2+l(t)\partial_{x_2}^2+b_1(t,x_1,x_2)\partial_{x_1}+b_2(t,x_1,x_2)\partial_{x_2}+c(t,x_1,x_2)\,.
$$
Moreover, the result can be extended to the operator $\mathcal{L}$ in (\ref{eq_L}) by considering the problem solved by $u(t,x_1,x_2)e^{-x_1^2-x_2^2}$, thus obtaining the following theorem.

\begin{teor}
There exist coefficients $a_{i,j}$, depending only on $t$, which are H\"older continuous of every order but not Lipschitz continuous and there exist $u\in\mathcal{H}_0$ such that the solution of problem (\ref{eq_probl_cauchy}) with $u_0=0$ and $f=0$ is not identically zero.$\hfill\square$
\end{teor}

In view of the previous results, a question naturally arises: which is the \emph{minimal} regularity with respect to $t$ (between Lipschitz continuity and H\"older continuity) of the coefficients of the principal part of $\mathcal{L}$ guaranteeing uniqueness of the solution of (\ref{eq_probl_cauchy})? To answer to this question, we recall the definition of \emph{modulus of continuity} that can be exploited to measure the degree of regularity of a function.

\begin{defi}
A \emph{modulus of continuity} is a function $\mu:[0,1]\to[0,1]$ which is continuous, increasing, concave and such that $\mu(0)=0$. A function $f:\mathbb{R}\to\mathbb{R}$ has \emph{regularity $\mu$} if
$$
\sup_{0<\vert t-s\vert<1}\frac{f(t)-f(s)}{\mu(t-s)}<+\infty\,.
$$
The set of all functions having regularity $\mu$ is denoted by $C^\mu$.
\end{defi}

As particular cases, the Lipschitz continuity, the $\tau$-H\"older continuity ($\tau\in]0,1[$) and the \emph{logarithmic Lipschitz} (in short \emph{Log-Lipschitz}) continuity are obtained for $\mu(s)=s$, $\mu(s)=s^{\tau}$ and $\mu(s)=s\log(1+1/s)$, respectively.

A further characterization of the modulus of continuity is the \emph{Osgood condition} which is crucial in most of the results on uniqueness and continuity that are described in the rest of the article. A modulus of continuity $\mu$ satisfies the Osgood condition if
$$
\int_0^1 \frac{1}{\mu(s)}ds = +\infty\,.
$$

This characterization is used, for instance, in~\cite{JMPA} to obtain the following result.

\begin{teor}\label{teo_delpri_2005}
Le $\mu$ be a modulus of continuity that satisfies the Osgood condition. Let
\begin{equation}\label{eq_H_1}
\mathcal{H}_1\triangleq H^1([0, T],L^2(\mathbb{R}^n))\cap L^2([0, T],H^2(\mathbb{R}^n))
\end{equation}
and let the coefficients $a_{i,j}$ in (\ref{eq_L}) be such that, for all $i,j=1,\ldots,n$,
$$
a_{i,j}\in C^\mu([0,T],\mathscr{C}_b(\mathbb{R}^n))\cap \mathscr{C}([0,T],\mathscr{C}_b^2(\mathbb{R}^n))\,,
$$
where $\mathscr{C}_b^2$ is the space of the bounded functions whose first and second derivatives are bounded. If $u\in \mathcal{H}_1$, if $\mathcal{L}u=0$ on $[0,T]\times \mathbb{R}^n$ and if $u(0,x)=0$ on $\mathbb{R}^n$, then $u\equiv 0$ on $[0,T]\times\mathbb{R}^n$.
\end{teor}

More recently, by using Bony's para-multiplication, the result has been improved as far as the regularity with respect to $x$ is concerned, i.e. replacing $\mathscr{C}^2$ regularity with Lipschitz regularity~\cite{AMPA}.

Note that the claim of Theorem~\ref{teo_delpri_2005} refers to the function space defined by (\ref{eq_H_1}), however, it is not difficult to extend it to the function space $\mathcal{H}_0$ defined by (\ref{eq_insieme}).

\section{Conditional stability results}
As mentioned in the introduction, for Cauchy problems related to the backward parabolic differential operators, which in general are not well posed, the notion of continuous dependence from initial data is replaced by the notion of (conditional) stability which is associated with the property of a problem to be well behaved, as defined by John~\cite{John}. The question about the conditional stability can be stated as follows. Suppose that two functions $u$ and $v$, defined in $[0,T]\times\mathbb{R}^n$, are solutions of the same equation; suppose, in addition, that $u$ and $v$ satisfy a fixed bound in a space $\mathcal{K}$ and that $\Vert u(0,\cdot)-v(0,\cdot)\Vert_{\mathcal{H}}$ is small (less than some $\epsilon$). Given these assumptions can we say something on the quantity $\sup_{t\in[0,T^\prime]}\Vert u(t,\cdot)-v(t,\cdot)\Vert_{\mathcal{K}}$ for some $T^\prime<T$? Does it remains small as well (e.g. less than a value related to $\epsilon$)? In this section we report some results that give an answer to the above questions.

\subsection{Stability with Lipschitz-continuous (with respect to $t$) coefficients}
One of the first results on conditional stability has been proven by Hurd~\cite{Hurd} in the same theoretical framework considered by Lions and Malgrange.

\begin{teor}\label{teor_LioMal}
Suppose that the coefficients $a_{i,j}$ are Lipschitz continuous in $t$ and in $x$. For every $T^\prime\in]0,T[$ and for every $D>0$ there exist $\rho>0$, $\delta\in]0,1[$ and $M>0$ such that if $u\in \mathcal{H}_0$ is a solution of $\mathcal{L}u=0$ on $[0,T]$ with $\Vert u(t,\cdot)\Vert_{L^2}<D$ on $[0,T]$ and $\Vert u(0,\cdot)\Vert_{L^2}<\rho$, then
\begin{equation}\label{eq_teo_hurd}
\sup_{t\in[0,T^\prime]}\Vert u(t,\cdot)\Vert_{L^2}\le M\Vert u(0,\cdot)\Vert_{L^2}^\delta\,.
\end{equation}
The constants $\rho$, $\delta$ and $M$ depend only on $T^\prime$ and $D$, on the ellipticity constant of $\mathcal{L}$, on the $L^\infty$ norms of the coefficients $a_{i,j}$, $b_j$, $c$ and of their spatial derivatives, and on the Lipschitz constant of the coefficients $a_{i,j}$ with respect to time.$\hfill\square$
\end{teor}

The result expressed by equation (\ref{eq_teo_hurd}) implies uniqueness of the solution to the Cauchy problem, so that a necessary condition to this kind of of conditional stability is that the coefficients $a_{i,j}$ fulfil the Osgood condition with respect to time. Hence a natural question arises: is Osgood condition also a sufficient condition? Del Santo and Prizzi~\cite{MATAN} have given a negative answer to this question. In particular, mimicking Pli\'s counterexample, they have shown that if the coefficients $a_{i,j}$ are not Lipschitz-continuous but only Log-Lipschitz-continuous then Hurd's result does not hold. Moreover, they have proven that is the coefficients are Log-Lipschitz-continuous then a conditional stability property, although weaker than (\ref{eq_teo_hurd}), does hold. More recently, the result has been further improved~\cite{NONLINAL}.

\subsection{Stability with Log-Lipschitz-continuous (with respect to $t$) coefficients}
As mentioned above, Osgood condition is not a sufficient condition for conditional stability of the solution. The following paragraph specifies this claim.
\subsubsection{Counterexample for the Lipschitz continuity case}\label{par_controesempio_uno}
The counterexample relies on the fact that it is possible~\cite{MATAN} to construct
\begin{itemize}
\item a sequence $\{\mathcal{L}_k\}_{k\in \mathbb{N}}$ of backward uniformly parabolic operators with uniformly Log-Lipschitz-continuous coefficients (not depending on the space variables) in the principal part and space-periodic uniformly bounded smooth coefficients in the lower order terms,
\item a sequence $\{u_k\}_{k\in\mathbb{N}}$ of space-periodic smooth uniformly bounded solutions of $\mathcal{L}_ku_k=0$ on $[0,1]\times\mathbb{R}^2$,
\item a sequence $\{t_k\}_{k\in\mathbb{N}}$ of real numbers, with $t_k\to 0$ as $k\to\infty$,
\end{itemize}
such that
$$
\lim_{k\to\infty}\Vert u_k(0,\cdot,\cdot)\Vert_{L^2([0,2\pi]\times[0,2\pi])}=0
$$
and
$$
\lim_{k\to\infty}\frac{\Vert u_k(t_k,\cdot,\cdot)\Vert_{L^2([0,2\pi]\times[0,2\pi])}}{\Vert u_k(0,\cdot,\cdot)\Vert^\delta_{L^2([0,2\pi]\times[0,2\pi])}}=+\infty
$$
for every $\delta>0$.
We remark that this situation is exactly what is needed to show that for backward operators with Log-Lipscitz continuous coefficient a result similar to Theorem~\ref{teor_LioMal} cannot hold.

\subsubsection{Stability result in the Log-Lipschitz case}
In the case of Log-Lipschitz coefficients a result weaker that (\ref{eq_teo_hurd}) is valid.

Consider the equation $\mathcal{L}u=0$ on $[0,T]\times\mathbb{R}^n$ and suppose that
\begin{enumerate}
\item for all $(t,x)\in[0,T]\times\mathbb{R}^n$ and for all $i,j=1,\ldots,n$, $a_{i,j}(t,x)=a_{j,i}(t,x)$;\label{condizione_prima}
\item there exists $k>0$ such that, for all $(t,x,ξ)\in[0,T]\times\mathbb{R}^n\times\mathbb{R}^n$,
$$
k\vert\xi\vert^2\le\sum_{i,j=1}^n a_{i,j}(t,x)\xi_i\xi_j\le k^{−1}\vert\xi\vert^2\,;
$$
\item for all $i,j=1,\ldots,n$, $a_{i,j}\in \textnormal{LogLip}([0,T],L^\infty(\mathbb{R}^n))\cap L^\infty([0,T],\mathscr{C}^2_b(\mathbb{R}^n))$, in particular
$$
\sup_{x\in\mathbb{R}^n,0<\vert\tau\vert<1}\frac{\vert a_{i,j}(t+\tau,x)-a_{i,j}(t,x)\vert}{\vert\tau\vert\left(\log\left(1+\frac{1}{\vert\tau\vert}\right)\right)}<+\infty\,;\label{condizione_terza}
$$
\item for all $j=1,\ldots,n$, $b_j\in L^\infty([0,T],\mathscr{C}_b^2(\mathbb{R}^n))$;\label{condizione_quarta}
\item $c\in L^\infty([0,T],\mathscr{C}^2_b(\mathbb{R}^n))$.\label{condizione_ultima}
\end{enumerate}

\begin{teor}\label{teo_matan}
{\bf\cite{MATAN}}
Suppose that the above hypotheses \ref{condizione_prima}-\ref{condizione_ultima} hold. For all $T^\prime\in]0,T[$ and for all $D>0$ there exist $\rho>0$, $M>0$, $N>0$ and $0<\beta<1$ such that, if $u\in\mathcal{H}_0$ is a solution of $\mathcal{L}u=0$ on $[0,T]$ with
$$
\sup_{t\in[0,T]}\Vert u(t,\cdot)\Vert_{L^2}\le D
$$
and $\Vert u(0,\cdot)\Vert_{L^2}\le \rho$, then
\begin{equation}\label{eq_teo_del_pri}
\sup_{t\in[0,T^\prime]}\Vert u(t,\cdot)\Vert_{L^2}\le M e^{-N\vert\log\Vert u(0,\cdot)\Vert_{L^2}\vert^\beta}\,,
\end{equation}
where the constants $\rho$, $\beta$, $M$ and $N$ depend only on $T^\prime$, on $D$, on the ellipticity constant of $\mathcal{L}$, on the $L^\infty$ norms of the coefficients $a_{i,j}$ and of their spatial first derivatives, and on the Log-Lipschitz constant of the coefficients $a_{i,j}$ with respect to time.
\end{teor}

Using Bony's para-product the result can be extended to the case in which the coefficient are not necessarily $\mathscr{C}_b^2$-continuous with respect to $x$ but only Lipschitz~\cite{NONLINAL}.


\subsection{Stability with Osgood-continuous (with respect to time) coefficients}
Let us finally come to the new result contained in this paper. As in the previous section we first present a counterexample to the stability condition (\ref{eq_teo_del_pri}) and then a new weaker stability result.

\subsubsection{Counterexample for the Log-Lipschitz case}\label{par_controesempio_due}
Consider the modulus of continuity $\omega$ defined by
$$
\omega(s)=s\log\left(1+\frac{1}{s}\right)\log\left(\log\left(1+\frac{1}{s}\right)\right)
$$
and note that $\omega$ satisfies the Osgood condition but is not Log-Lipschitz continuous. Analogously to Paragraph~\ref{par_controesempio_uno}, it is possible to construct (see~\cite{mia_tesi})
\begin{itemize}
\item a sequence $\{\mathcal{P}_k\}_{k\in \mathbb{N}}$ of backward uniformly parabolic operators with uniformly $\mathscr{C}^\omega$-continuous coefficients (not depending on the space variables) in the principal part and space-periodic uniformly bounded smooth coefficients in the lower order terms,
\item a sequence $\{u_k\}_{k\in\mathbb{N}}$ of space-periodic smooth uniformly bounded solutions of $\mathcal{L}_ku_k=0$ on $[0,1]\times\mathbb{R}^2$,
\item a sequence $\{t_k\}_{k\in\mathbb{N}}$ of real numbers, with $t_k\to 0$ as $k\to\infty$,
\end{itemize}
such that
$$
\lim_{k\to\infty}\Vert u_k(0,\cdot,\cdot)\Vert_{L^2([0,2\pi]\times[0,2\pi])}=0
$$
but (\ref{eq_teo_del_pri}) does not hold for all $k$; more precisely
$$
\lim_{k\to\infty}\frac{\Vert u_k(t_k,\cdot,\cdot)\Vert_{L^2([0,2\pi]\times[0,2\pi])}}{e^{-N\vert \log\Vert u_k(0,\cdot,\cdot)\Vert_{L^2([0,2\pi]\times[0,2\pi])}\vert^\delta}}=+\infty
$$
for every $\delta>0$.

\subsubsection{Stability result in the Osgood-continuous case}

Let $\mathcal{L}$ be a backward parabolic operator whose coefficients depend only on $t$, i.e. let
$$
\mathcal{L}u=\partial_t u+\sum_{i,j=1}^na_{i,j}(t)\partial_{x_i}\partial_{x_j}u+\sum_{j=1}^nb_j(t)\partial_{x_j}u+c(t)u
$$
on the strip $[0,T]\times\mathbb{R}^n$. Suppose that $a_{i,j}(t)=a_{j,i}(t)$ for all $i,j=1,\ldots,n$ and for all $t\in[0,T]$. Let $a_{i,j},b_j,c\in L^\infty([0,T])$, for all $i,j=1,\ldots,n$. Let $\mu$ be a modulus of continuity satisfying the Osgood condition. Let $a_{i,j}$ $\mathscr{C}^\mu$-continuous, i.e.
$$
\sup_{0<\vert\tau\vert<1}\frac{\vert a_{i,j}(t+\tau)-a_{i,j}(t)\vert}{\mu(\vert\tau\vert)}<+\infty\,.
$$

\begin{teor}\label{teo_nuovo}
For all $T^\prime\in]0,T[$ and for all $D>0$ there exist $\rho>0$, and there exists an increasing continuous function $\Psi:[0,+\infty)\to[0,+\infty)$, with $\Psi(0)=0$ such that, if $u\in\mathcal{H}_0$ is a solution of $\mathcal{L}u=0$ on $[0,T]$ with $\Vert u(t,\cdot)\Vert_{L^2}\le D$ on $[0,T]$ and $\Vert u(0,\cdot)\Vert_{L^2}\le\rho$, then
\begin{equation}\label{eq_in_L2}
\sup_{t\in[0,T^\prime]}\Vert u(t,\cdot)\Vert_{L^2}\le\Psi(\Vert u(0,\cdot)\Vert_{L^2})\,.
\end{equation}
The constant $\rho$ and the function $\Psi$ depend only on $T^\prime$, on $D$, on the ellipticity constant of $\mathcal{L}$, on the $L^\infty$ norms of the coefficients $a_{i,j}$ and of their spatial first derivatives, and on the Osgood constant of the coefficients $a_{i,j}$.$\hfill\square$
\end{teor}

\subsubsection{Comments on the result and its proof}
The complete proof of Theorem~\ref{teo_nuovo} is beyond the aims of this paper and is not reported here. However, to provide the reader with some insights about the demonstration, in the following we comment on the analogies and the differences between the new result and the previous ones~\cite{NONLINAL,MATAN}. We begin by recalling that Theorem~\ref{teo_matan} 
is a consequence of the ``energy'' estimate (see Proposition 1 in~\cite{MATAN})
\begin{multline}\label{eq_stima_vecchia}
\int_0^s e^{2\gamma t}e^{-2\beta\phi_\lambda((t+\tau)/\beta)}\Vert u(t,\cdot)\Vert^2_{H^{1-\alpha t}}dt\le\\[2mm]
\le M\big((s+\tau)e^{2\gamma s}e^{2\beta\phi_\lambda((s+\tau)/\beta)}\Vert u(s,\cdot)\Vert^2_{H^{1-\alpha s}}+\\[2mm]
+\tau\phi_\lambda^\prime(\tau/\beta)e^{-2\beta\phi_\lambda(\tau/\beta)}\Vert u(0,\cdot)\Vert^2_{L^2}\big)\,,
\end{multline}
where $\phi_\lambda$ is the solution of the differential equation
\begin{equation}\label{eq_phi_vecchia}
y\phi_\lambda^{\prime\prime}(y)=-\lambda\phi_\lambda^\prime(y)(1+\vert \log\phi^\prime(y)\vert)\,.
\end{equation}
and the constants depend, in particular\footnote{They also depend, as specified in the claim of the theorem, on $T^\prime$, on $D$, on the ellipticity constant of $\mathcal{L}$, on the $L^\infty$ norms of the coefficients $a_{i,j}$ and of their spatial first derivatives. The parameter $\lambda$ also depends on these quantities.}, on the Log-Lipschitz constant of the coefficients $a_{i,j}$ with respect to time. Now, the novelty of Theorem~\ref{teo_nuovo} is that the coefficients $a_{i,j}$ are supposed to be only Osgood-continuous, hence there is no Log-Lipschitz constant to be taken as a reference. On the other hand, the energy estimate will necessarily contain information on the modulus of continuity (which is assumed to verify the Osgood condition). Indeed, the energy estimate is
\begin{multline}\label{eq_u_hat_seconda}
\frac{1}{4}\left(k_A\vert\xi\vert^2+\gamma\right)\int_0^\sigma e^{(1-\alpha t)\vert\xi\vert^2\omega\left(\frac{1}{\vert\xi\vert^2+1}\right)}e^{2\gamma t}e^{-2\beta\phi_\lambda\left(\frac{t+\tau}{\beta}\right)}\vert\hat{u}(t,\xi)\vert^2dt\le\\
\le \phi^\prime_\lambda\left(\frac{\tau}{\beta}\right)\tau e^{\vert\xi\vert^2\omega\left(\frac{1}{\vert\xi\vert^2+1}\right)}e^{-2\beta\phi_\lambda\left(\frac{\tau}{\beta}\right)}\vert\hat{u}(0,\xi)\vert^2+\\
+(\sigma+\tau)(\gamma+k_A^{-1}\vert\xi\vert^2)e^{2\gamma\sigma}e^{-2\beta\phi_\lambda\left(\frac{\sigma+\tau}{\beta}\right)}\vert\hat{u}(\sigma,\xi)\vert^2\,,
\end{multline}
where, in particular, $\hat{u}$ denotes the Fourier transform of $u$ with respect to $x$, $\omega$ is the modulus of continuity of the coefficients $a_{i,j}$, $k_A$ is the ellipticity constant of the principal part of $\mathcal{L}$ and $\phi_\lambda$ is now the solution of the differential equation
\begin{equation}\label{eq_phi_nuova}
y\phi^{\prime\prime}_\lambda(y)=-\lambda(\phi^{\prime}_\lambda(y))^2\omega\left(\frac{k_A}{\phi^\prime_\lambda(y)}\right)\,,
\end{equation}
where, again, the modulus of continuity appears. By comparing (\ref{eq_stima_vecchia})-(\ref{eq_phi_vecchia}) with (\ref{eq_u_hat_seconda})-(\ref{eq_phi_nuova}) one can see that Theorem~\ref{teo_nuovo} is not a trivial generalization of Theorem~\ref{teo_matan}. In addition, (\ref{eq_u_hat_seconda}) leads by integrating in $\xi$ to the estimate
\begin{equation}\label{eq_stima_prima}
\sup_{z\in[0,\bar{\sigma}]}\Vert u(z,\cdot)\Vert_{H_{\frac{1}{2},\omega}^1}^2\le Ce^{-\sigma\phi_\lambda^\prime\left(\frac{\sigma+\tau}{\beta}\right)}\left[\phi^\prime_\lambda\left(\frac{\tau}{\beta}\right)e^{-2\beta\phi_\lambda\left(\frac{\tau}{\beta}\right)}\Vert u(0,\cdot)\Vert^2_{H^0_{1,\omega}}+\Vert u(\sigma,\cdot)\Vert_{H^1}\right]\,,
\end{equation}
where, in particular, the function spaces $H^1_{\frac{1}{2},\omega}$ and $H^0_{1,\omega}$ come into the scene. These spaces, defined by
$$
\Vert u\Vert^2_{H^d_{a,\omega}}\triangleq \int_{\mathbb{R}^n}\left(1+\vert \xi\vert^2\right)^de^{a\vert\xi\vert^2\omega\left(\frac{1}{\vert\xi\vert^2+1}\right)}\vert \hat{u}(\xi)\vert^2d\xi<+\infty\,,
$$
are tailored on the modulus of continuity $\omega$ and, although comparable with Gevrey-Sobolev spaces, they do not coincide with any of them.

The final estimate (\ref{eq_in_L2}), which is written with respect to the $L^2$ norm, can be obtained from (\ref{eq_stima_prima}) by exploiting the regularising properties of the (forward) parabolic operator.

\bibliographystyle{plain}
\bibliography{ISAAC_bibliografia}

\end{document}